\documentclass[11pt]{article}
\usepackage{fullpage}
\usepackage{amsmath}
\usepackage{amssymb}
\usepackage{dsfont}
\usepackage[final]{graphicx}
\usepackage{graphicx}
\usepackage{amsthm}
\usepackage{color}
\usepackage{bbm}
\usepackage{setspace}
\usepackage{multirow}
\usepackage{ifpdf}
\usepackage{float}
\usepackage[small]{caption}
\usepackage[pdftex]{thumbpdf}
\usepackage{appendix}
\usepackage[pdftex,colorlinks]{hyperref}
\hypersetup{allcolors=blue}
 \usepackage{natbib}
\restylefloat{figure}
\theoremstyle{plain}              

\numberwithin{definition}{section}
\theoremstyle{remark} 

\numberwithin{equation}{section}

\newcommand{\e}{\mathbb{E}}
\newcommand{\p}{\mathbb{P}}

\numberwithin{equation}{section}
\date{}
\begin{document}
\title{A 3-Queue Polling System with Join the Shortest - Serve the Longest Policy}
\author{Efrat Perel$^{1}$, Nir Perel$^{1}$\footnote{Corresponding author, e-mail: nirp@afeka.ac.il} \ and  Uri Yechiali$^{2}$}
\maketitle
\begin{center}
\small $^1$School of Industrial Engineering and Management, \\
Afeka Tel-Aviv Academic College of Engineering,
Tel-Aviv, Israel. \\

\small $^2$Department of Statistics and Operations Research,\\
School of Mathematical Sciences,
Tel-Aviv University, Tel-Aviv, Israel.
\end{center}

\begin{abstract}
In 1987, J.W. Cohen analyzed the so-called Serve the Longest Queue (SLQ) queueing system, where a single server attends two non-symmetric $M/G/1$-type queues, exercising a non-preemptive priority switching policy. Cohen further analyzed in 1998 a non-symmetric 2-queue Markovian system, where newly arriving customers follow the Join the Shortest Queue (JSQ) discipline. The current paper generalizes and extends Cohen's works by studying a combined JSQ-SLQ model, and by broadening the scope of analysis to a non-symmetric 3-queue system, where arriving customers follow the JSQ strategy and a single server exercises the preemptive priority SLQ discipline. The system states' multi-dimensional probability distribution function is derived while applying a non-conventional representation of the underlying process's state-space. The analysis combines both Probability Generating Functions and Matrix Geometric methodologies. It is shown that the joint JSQ-SLQ operating policy achieves extremely well the goal of balancing between queue sizes. This is emphasized when calculating the Gini Index associated with the differences between mean queue sizes: the value of the coefficient is close to zero. Extensive numerical results are presented.
\end{abstract}


\section{Introduction}\label{SL3 intro}
Queueing systems with multiple queues have been studied extensively in the literature. Noted are models where, upon arrival, customers join the shortest queue (JSQ), along with other models where the server serves the longest queue (SLQ). The main goal of each model-type is to balance between the queue sizes: in the JSQ model, each customer's action is directed into achieving the balance, while in the SLQ model, the server's actions are aimed at this goal.

\cite{cohen1987two} studied a single server system with two $M/G/1$-type queues, under a non-preemptive Serve the Longest Queue (SLQ) regime. Each queue has its own generally distributed service time and its own Poisson rate of arrival, where upon each service completion, the server selects the longer queue. If both queues are equal in size, the server chooses the next queue to be served according to some probability. This model is analyzed by constructing and solving a boundary value problem.
The SLQ regime was further studied by \cite{flatto1989longer} who analyzed a Markovian system with two identical queues and server's preemptive switching policy, by calculating a 2-dimensional probability generating function of the queue lengths and studying its asymptotic behavior. \cite{houtum1997symmetric} analyzed a non-preemptive Markovian system with $N$ symmetric queues, assuming that the difference between the longest and the shortest queue is limited by some known constant. As a result, it became possible to perform exact analysis of the model and its variants via matrix geometric methods.  \cite{knessl2013nonsymmetric} derived asymptotic properties for a Markovian non-symmetric two-queue model in heavy traffic, in which the server always attends the longer queue, and in case when both queues are with equal sizes, the server splits its capacity between the queues. An exact integral representation is derived for the joint probability distribution of the queue sizes.  \cite{ravid2013repair} studied a single repairman problem with two Markovian queues having non-identical arrival rates, and pointed their model's relation to a corresponding SLQ model. The SLQ model was also studied in the context of wireless communications and switching scheduling. For example, \cite{maguluri2014stability} analyzed a wireless network with SLQ mechanism, while \cite{pedarsani2016stability} investigated the stability of an open multi-class network under SLQ modeling.

The Join the Shortest Queue (JSQ) model has been also investigated widely in the literature. \cite{winston1977optimality} showed that in a symmetric Markovian system with identical servers, the JSQ policy maximizes the discounted number of jobs completing service by some time $t$.
\cite{halfin1985shortest} considered a service center with two parallel single-server queues and identical Exponential service rates, and derived, using linear programming techniques, bounds for the probability distribution of the number of customers in the system. \cite{hordijk1990optimality} showed, using Markov Decision Processes, that the JSQ policy is optimal for queues with finite buffers and a general arrival process. \cite{adan1991analysis} studied a system with two non-symmetric parallel queues, and showed, by using an elementary compensation procedure, that the equilibrium distribution of the lengths of the two queues can be exactly represented by an infinite sum of product form solutions. In a following paper \cite{adan1991analysisjockeying} considered the same system with jockeying between the queues, and analyzed the model using a matrix geometric approach. \cite{cohen1998analysis} presented an analytic solution for an asymmetric two-server queueing model with Poisson arrivals, where arriving customers join the shorter queue and service times are Exponentially distributed, by determining explicitly the bivariate generating function of the stationary joint distribution of the queue lengths. \cite{foley2001join} derived asymptotic results for a Markovian system with two servers, each of which receives a dedicated stream of customers in addition to customers from a stream of smart customers who join the shorter queue. \cite{yao2005infinite,yao2006infinite} provided an asymptotic analysis for a model with two $M/M/\infty$ queues where a new arrival is routed to the queue with the smaller number of customers, by establishing and analyzing a bi-variate generating function for the number of customers in each queue. \cite{blanc2009bad} studied a non-symmetric multi-server Markovian system where customers follow the JSQ policy, and quantified the probability that a customers's service would have started earlier if one had joined another queue than the queue that was actually chosen. By using the compensation procedure, \cite{adan2013erlang} analyzed a system with two single-server queues, JSQ policy, where customers inter-arrival times follow an Erlang distribution and service times are Exponentially distributed. \cite{adan2016shorter} analyzed the JSQ polling model by using the compensation method and by solving a boundary value problem. \cite{eschenfeldt2018join} and \cite{braverman2020steady}  analyzed the steady-state properties of the JSQ model under the Halfin-Whitt regime.  \cite{dawson2019performance} studied a queueing network with a large number of nodes, in which each queue has a dedicated input stream, while, in addition, there is an extra input stream, balancing the network load by directing its arrivals to the shortest queue. A mean field interaction model is set up to study the performance of this network in terms of limit results. \cite{dimitriou2021analysis} studied a communication-related Markovian JSQ model with two infinite capacity orbit queues. An arriving job that finds the server busy is forwarded to the smaller-sized orbit queue. It is worth noting that in most of the JSQ models described above, \emph{each queue has its own server}. In contrast, in this work, the combined JSQ-SLQ model analyzes a single-server polling-type system, where the server switches from a shorter queue to a longer one, while arriving customers join the shortest queue.

Recently, \cite{perel2020polling} presented a combined unifying Join the shortest queue - Serve the longest Queue (JSQ-SLQ) polling-type model under a preemptive priority switching policy of the server. Such a combined policy aims at equalizing the queue sizes. Specifically, the authors considered a polling system comprised of two non-identical Markovian queues, attended by a single server that alternates between them. An arriving customer always joins the shortest queue, unless the queue lengths are equal, in which case the customer joins queue $i$ with probability $p_i\geq0$, $i=1,2$, $p_1+p_2=1$. On the other hand, the server always serves the longest queue while exercising a queue-size depending preemptive priority policy, i.e. the server never resides in a shorter queue, giving priority to the longest queue. This preemptive-type policy implies that at a moment when the number of customers in an un-served queue exceeds the number of customers in the served queue (either at a service completion or when arrival occurs), the server immediately switches to the longest queue.

The current paper also combines the above two operating policies, i.e. JSQ and SLQ, into a unified model, and extends the analysis to a 3-queue system. Specifically, we study a polling system comprised of three asymmetrical Markovian queues, denoted by $Q_1$, $Q_2$ and $Q_3$, and a single server that alternates between them according to a preemptive SLQ policy. That is, the server always attends the longest queue. As soon as the number of customers in one of the non-attended queues rises above the number of customers in the attended queue, the server stops serving the served customer and instantaneously switches to the longest queue. The service of the interrupted customer will resume anew when its turn comes again (within each queue the service order is FCFS). In addition, if upon service completion the two non-attendant queues ($Q_i$ and $Q_j$) have the (same) longest queue length, the server will switch to $Q_i$ $(Q_j)$ w.p. $\frac{q_i}{q_i+q_j}$ $\left(\frac{q_j}{q_i+q_j}\right)$, where $q_i>0$, $i=1,2,3$, is the weight the server assigns to $Q_i$. It is assumed that service duration of an arbitrary customer in $Q_i$ is Exponentially distributed with mean $1/\mu_i$, $i=1,2,3$. Customers arrive at the system according to a Poisson process with rate $\lambda$. Each arriving customer exercises the Join the Shortest Queue (JSQ) policy, whereas if the lengths of all queues are equal, the customer joins $Q_i$ w.p. $p_i$, where $p_1+p_2+p_3=1$. Furthermore, if only two queues, $Q_i$ and $Q_j$ ($i\neq j$), have the same shortest length, an arriving customer, being consistent, joins $Q_i$ $(Q_j)$ w.p. $\frac{p_i}{p_i+p_j}$ $\left(\frac{p_j}{p_i+p_j}\right)$ (note that according to this setting, at most one of the $p_i's$ can be zero).

Even though we consider here a 3-queue system, where each queue is un-bounded in its size, we are able to formulate this model as a 2-dimensional continuous-time Markovian process with only one un-bounded size. In many cases, solving a 2-dimensional un-bounded system may be done by applying a boundary value problem analysis (see e.g. \cite{flatto1989longer}, \cite{avrachenkov2014retrial}, \cite{adan2016shorter}), or by a truncation method as carried out in \cite{bright1995calculating}. In contrast, using an un-conventional approach, as was introduced in \cite{perel2020polling}, we are able to derive the equilibrium joint probability mass function of the three queue lengths. For that purpose, we apply both the two commonly used methods to solve QBD problems, namely, the probability generating functions method and the matrix geometric analysis. The main innovation of our analytical approach is to define only a finitely-many number of probability generating functions representing the system's states.

The paper is organized as follows. In Section \ref{SL3 sec: model} the model is described and formulated. In Section \ref{SL3 sec: pgf} the steady-state equations are established and the corresponding probability generating functions are derived. Various performance measures are defined and calculated. In Section \ref{SL3 sec: MG} we utilize the matrix geometric method in order to derive the system's stability condition, as well as the equilibrium distribution of the system's states. Extensive numerical results are presented in Section \ref{SL3 sec results}, indicating that the balance between the 3 queue sizes is very well achieved. This is further emphasized by the corresponding Gini index that its value is close to zero in almost all cases examined. Section \ref{SL3 sec conclusions} concludes the paper.

\section{The Model}\label{SL3 sec: model}
We consider a polling system with a single server and \textmd{three} non-identical queues, denoted by $Q_1$, $Q_2$ and $Q_3$. Customers arrive at the system according to a Poisson process with rate $\lambda$. Each arriving customer exercises the 'Join the Shortest Queue' (JSQ) policy, whereas, if the lengths of all queues are equal, the customer joins $Q_i$ w.p. $p_i\geq 0$, where $p_1+p_2+p_3=1$. Furthermore, if only two queues, $Q_i$ and $Q_j$ ($i\neq j$), have the same shortest length, an arriving customer, being consistent, joins $Q_i$ $(Q_j)$ w.p. $\frac{p_i}{p_i+p_j}$ $\left(\frac{p_j}{p_i+p_j}\right)$ (note that according to this setting, at most one of the $p_i's$ can be zero). Service duration of an arbitrary customer in $Q_i$ is Exponentially distributed with mean $1/\mu_i$, $i=1,2,3$. The server alternates between the three queues according to the preemptive 'Serve the Longest Queue' (SLQ) switching policy. That is, the server always attends (one of) the longest queue(s). As soon as the number of customers in one of the non-attendant queues rises above the number of customers in the attendant queue, the server stops serving the served customer and instantaneously switches to the longest queue. The service of the interrupted customer will resume anew when its turn comes again (within each queue the service order is FCFS). In addition, if upon service completion the two non-attendant queues ($Q_i$ and $Q_j$) have the (same) longest queue length, the server will switch to $Q_i$ $(Q_j)$ w.p. $\frac{q_i}{q_i+q_j}$ $\left(\frac{q_j}{q_i+q_j}\right)$, where $q_i>0$, $i=1,2,3$, is the weight assigned to $Q_i$. If, upon service completion at $Q_i$ all queues are equal in size, the server stays in $Q_i$.

At time $t>0$, let $L_i(t)$ denote the number of customers present in $Q_i$, $i=1,2,3$, and, assuming stability (see in the sequel analysis of the stability condition), let $L_i=\lim\limits_{t\rightarrow\infty}L_i(t)$. The conventional approach would be to define a 3-dimensional un-bounded space of the system's states $\{L_1(t),L_2(t),L_3(t)\}$ which seems to be untractable. Instead, we construct a non-conventional state space (following \cite{perel2020polling}) as follows: Let $I(t)$ denote the index of the queue attended by the server at time $t$, and define the triple $D(t)=\left(I(t);D_2(t),D_3(t)\right)$, where $D_j(t)=L_1(t)-L_j(t)$, $j=2,3$. Let $D=\lim\limits_{t\rightarrow\infty}D(t)$.

We now formulate the 3-dimensional intriguing polling system as a two-dimensional continuous-time Markovian process, with state space $\{(n,d)\}$, where $n\in\{0,1,2,3,...\}$, and
\begin{align*}
d\in\mathfrak{D}=\{&(1;1,1),(1;0,1),(1;1,0),(1;0,0),(2;0,1),(2;0,0),(2;-1,0),(2;-1,-1),\\
                  &(3;1,0),(3;0,0),(3;0,-1),(3;-1,-1)\}.
\end{align*}
In steady state, the joint probability distribution function is denoted by $P_{n,d}=\p(L_1=n,D=d)$.

\section{Steady-state analysis using probability generating functions}\label{SL3 sec: pgf}
In this section we construct the steady-state joint probability distribution function of the process $\{(n,d)\}$. We use the above mentioned non-conventional construction (see  \cite{perel2020polling}) to obtain $12$ partial probability generating functions (PGFs) defining the system's states and utilize their properties, as described below.
\subsection{Balance equations and PGFs}\label{SL3 sec balance equations}
Writing the balance equations for all $n$, we obtain:\\
For $d=(1;1,1)$,
\begin{align}
(\lambda+\mu_1) P_{n,(1;1,1)}=&\lambda p_1(P_{n-1,(1;0,0)}+P_{n-1,(2;0,0)}+P_{n-1,(3;0,0)})+\mu_2P_{n,(2;0,1)}+\mu_3P_{n,(3;1,0)},  &n\geq1. \label{SL3 eq: d=(1;1,1)}
\end{align}
For $d=(1;0,1)$,
\begin{align}
&(\lambda+\mu_1)P_{n,(1;0,1)}=\lambda\frac{p_2}{p_2+p_3}P_{n,(1;1,1)}+\mu_3\frac{q_1}{q_1+q_2}P_{n,(3;0,0)},  &n\geq1. \label{SL3 eq: d=(1;0,1)}
\end{align}
For $d=(1;1,0)$,
\begin{align}
&(\lambda+\mu_1)P_{n,(1;1,0)}=\lambda\frac{p_3}{p_2+p_3}P_{n,(1;1,1)}+\mu_2\frac{q_1}{q_1+q_3}P_{n,(2;0,0)},  &n\geq1.  \label{SL3 eq: d=(1;1,0)}
\end{align}
For $d=(1;0,0)$,
\begin{align}
&\lambda P_{0,(1;0,0)}=\mu_1P_{1,(1;1,1)}, \ &n=0. \label{SL3 eq: d=(1;0,0)}\\
&(\lambda+\mu_1)P_{n,(1;0,0)}=\lambda \left(P_{n,(1;1,0)}+P_{n,(1;0,1)}\right)+\mu_1P_{n+1,(1;1,1)},  &n\geq1. \label{SL3 eq: d(1;0,0) n}
\end{align}

For $d=(2;0,1)$,
\begin{align}
&(\lambda+\mu_2) P_{n,(2;0,1)}= \lambda \frac{p_1}{p_1+p_3}P_{n-1,(2;-1,0)}+\mu_3\frac{q_2}{q_1+q_2}P_{n,(3;0,0)},  &n\geq1. \label{SL3 eq: d=(2;0,1)}
\end{align}
For $d=(2;0,0)$,
\begin{align}
&\lambda P_{0,(2;0,0)}=\mu_2P_{0,(2;-1,0)}, \ &n=0. \label{SL3 eq: d=(2;0,0)}\\
&(\lambda+\mu_2)P_{n,(2;0,0)}=\lambda \left(P_{n-1,(2;-1,-1)}+P_{n,(2;0,1)}\right)+\mu_2P_{n,(2;-1,0)},  &n\geq1. \label{SL3 eq: d(2;0,0) n}
\end{align}
For $d=(2;-1,0)$,
\begin{align}
(\lambda+\mu_2)P_{n,(2;-1,0)}=&\lambda p_2(P_{n,(1;0,0)}+P_{n,(2;0,0)}+P_{n,(3;0,0)})+\mu_1P_{n+1,(1;0,1)}+\mu_3P_{n,(3;-1,-1)},  &n\geq0. \label{SL3 eq: d=(2;-1,0)}
\end{align}
for $d=(2;-1,-1)$,
\begin{align}
&(\lambda+\mu_2)P_{n,(2;-1,-1)}=\lambda \frac{p_3}{p_1+p_3}P_{n,(2;-1,0)}+\mu_1\frac{q_2}{q_2+q_3}P_{n+1,(1;0,0)},  &n\geq0. \label{SL3 eq: d=(2;-1,-1)}
\end{align}

For $d=(3;1,0)$,
\begin{align}
&(\lambda+\mu_3) P_{n,(3;1,0)}= \lambda \frac{p_1}{p_1+p_2}P_{n-1,(3;0,-1)}+\mu_2\frac{q_3}{q_1+q_3}P_{n,(2;0,0)},  &n\geq1. \label{SL3 eq: d=(3;1,0)}
\end{align}
For $d=(3;0,0)$,
\begin{align}
&\lambda P_{0,(3;0,0)}=\mu_3P_{0,(3;0,-1)},  &n=0. \label{SL3 eq: d=(3;0,0)}\\
&(\lambda+\mu_3)P_{n,(3;0,0)}=\lambda \left(P_{n-1,(3;-1,-1)}+P_{n,(3;1,0)}\right)+\mu_3P_{n,(3;0,-1)},  &n\geq1. \label{SL3 eq: d(3;0,0) n}
\end{align}
For $d=(3;0,-1)$,
\begin{align}
&(\lambda+\mu_3)P_{n,(3;0,-1)}=\lambda p_3(P_{n,(1;0,0)}+P_{n,(2;0,0)}+P_{n,(3;0,0)})+\mu_1P_{n+1,(1;1,0)}+\mu_2P_{n,(2;-1,-1)},  &n\geq0. \label{SL3 eq: d=(3;0,-1)}
\end{align}
Finally, for $d=(3;-1,-1)$,
\begin{align}
&(\lambda+\mu_3)P_{n,(3;-1,-1)}=\lambda \frac{p_2}{p_1+p_2}P_{n,(3;0,-1)}+\mu_1\frac{q_3}{q_2+q_3}P_{n+1,(1;0,0)},  &n\geq0. \label{SL3 eq: d=(3;-1,-1)}
\end{align}
Note that in equations (\ref{SL3 eq: d=(1;1,1)})--(\ref{SL3 eq: d=(3;-1,-1)}) above, equations (\ref{SL3 eq: d=(1;0,0)}), (\ref{SL3 eq: d=(2;0,0)}) and (\ref{SL3 eq: d=(3;0,0)}) hold for $n=0$; equations (\ref{SL3 eq: d=(2;-1,0)}), (\ref{SL3 eq: d=(2;-1,-1)}), (\ref{SL3 eq: d=(3;0,-1)}) and (\ref{SL3 eq: d=(3;-1,-1)}) hold for $n\geq0$, while the rest of the equations hold for $n\geq1$.

For each $d\in\mathfrak{D}$, define the conditional 
probability generating function of the number of customers in $Q_1$ as: \[G_d(z)=\sum_{n}{P_{n,d}z^n}, \ \ \ d\in\mathfrak{D}.\]

Multiplying equation (\ref{SL3 eq: d=(1;1,1)}) by $z^n$ and summing over $n\geq1$ results in
\begin{align}\label{SL3 eq: G(1;1,1)(z)}
(\lambda+\mu_1)G_{(1;1,1)}(z)=&\lambda p_1zG_{(1;0,0)}(z)+\lambda p_1zG_{(2;0,0)}(z)+\lambda p_1zG_{(3;0,0)}(z)\nonumber\\
&+\mu_2G_{(2;0,1)}(z)+\mu_3G_{(3;1,0)}(z).
\end{align}

Repeating this process for all $d\in\mathfrak{D}$, while using equations (\ref{SL3 eq: d=(1;0,1)})--(\ref{SL3 eq: d=(3;-1,-1)}), leads to
\begin{align}
&(\lambda+\mu_1)G_{(1;0,1)}(z)=\lambda\frac{p_2}{p_2+p_3}G_{(1;1,1)}(z)+\mu_3\frac{q_1}{q_1+q_2}G_{(3;0,0)}(z)-\mu_3\frac{q_1}{q_1+q_2}P_{0,(3;0,0)}, \label{SL3 eq: G(1;0,1)(z)}
\end{align}
\begin{align}
&(\lambda+\mu_1)G_{(1;1,0)}(z)=\lambda\frac{p_3}{p_2+p_3}G_{(1;1,1)}(z)+\mu_2\frac{q_1}{q_1+q_3}G_{(2;0,0)}(z)-\mu_2\frac{q_1}{q_1+q_3}P_{0,(2;0,0)},\label{SL3 eq: G(1;1,0)(z)}
\end{align}
\begin{align}
&(\lambda+\mu_1)zG_{(1;0,0)}(z)=\lambda zG_{(1;1,0)}(z)+\lambda zG_{(1;0,1)}(z)+\mu_1G_{(1;1,1)}(z)+\mu_1zP_{0,(1;0,0)},\label{SL3 eq: G(1;0,0)(z)}
\end{align}
\begin{align}
&(\lambda+\mu_2)G_{(2;0,1)}(z)=\lambda z\frac{p_1}{p_1+p_3}G_{(2;-1,0)}(z)+\mu_3\frac{q_2}{q_1+q_2}G_{(3;0,0)}(z)-\mu_3\frac{q_2}{q_1+q_2}P_{0,(3;0,0)}, \label{SL3 eq: G(2;0,1)(z)}
\end{align}
\begin{align}
&(\lambda+\mu_2)G_{(2;0,0)}(z)=\lambda zG_{(2;-1,-1)}(z)+\lambda G_{(2;0,1)}(z)+\mu_2G_{(2;-1,0)}(z)+\mu_2P_{0,(2;0,0)},\label{SL3 eq: G(2;0,0)(z)}
\end{align}
\begin{align}
&(\lambda+\mu_2)zG_{(2;-1,0)}(z)=\lambda p_2zG_{(1;0,0)}(z)+\lambda p_2zG_{(2;0,0)}(z)+\lambda p_2zG_{(3;0,0)}(z)\nonumber\\
&\ \ \ \ \ \ \ \ \ \ \ \ \ \ \ \ \ \ \ \ \ \ \ \ \ \ \ \ \ \ \ +\mu_1G_{(1;0,1)}(z)+\mu_3zG_{(3;-1,-1)}(z), \label{SL3 eq: G(2;-1,0)(z)}
\end{align}
\begin{align}
&(\lambda+\mu_2)zG_{(2;-1,-1)}(z)=\lambda z\frac{p_3}{p_1+p_3}G_{(2;-1,0)}(z)+\mu_1\frac{q_2}{q_2+q_3}G_{(1;0,0)}(z)-\mu_1\frac{q_2}{q_2+q_3}P_{0,(1;0,0)}, \label{SL3 eq: G(2;-1,-1)(z)}
\end{align}
\begin{align}
&(\lambda+\mu_3)G_{(3;1,0)}(z)=\lambda z\frac{p_1}{p_1+p_2}G_{(3;0,-1)}(z)+\mu_2\frac{q_3}{q_1+q_3}G_{(2;0,0)}(z)-\mu_2\frac{q_3}{q_1+q_3}P_{0,(2;0,0)}, \label{SL3 eq: G(3;1,0)(z)}
\end{align}
\begin{align}
&(\lambda+\mu_3)G_{(3;0,0)}(z)=\lambda zG_{(3;-1,-1)}(z)+\lambda G_{(3;1,0)}(z)+\mu_3G_{(3;0,-1)}(z)+\mu_3P_{0,(3;0,0)},\label{SL3 eq: G(3;0,0)(z)}
\end{align}
\begin{align}
&(\lambda+\mu_3)zG_{(3;0,-1)}(z)=\lambda p_3zG_{(1;0,0)}(z)+\lambda p_3zG_{(2;0,0)}(z)+\lambda p_3zG_{(3;0,0)}(z)\nonumber\\
&\ \ \ \ \ \ \ \ \ \ \ \ \ \ \ \ \ \ \ \ \ \ \ \ \ \ \ \ \ \ \ +\mu_1G_{(1;1,0)}(z)+\mu_2zG_{(2;-1,-1)}(z), \label{SL3 eq: G(3;0,-1)(z)}
\end{align}
\begin{align}
&(\lambda+\mu_3)zG_{(3;-1,-1)}(z)=\lambda z\frac{p_2}{p_1+p_2}G_{(3;0,-1)}(z)+\mu_1\frac{q_3}{q_2+q_3}G_{(1;0,0)}(z)-\mu_1\frac{q_3}{q_2+q_3}P_{0,(1;0,0)}, \label{SL3 eq: G(3;-1,-1)(z)}
\end{align}

The set of $12$ equations (\ref{SL3 eq: G(1;1,1)(z)})-(\ref{SL3 eq: G(3;-1,-1)(z)}) can be written in a matrix form as
\begin{equation}\label{SL3 eq: set of PGF}
A(z)\cdot\vec{G}(z)=\vec{P}(z),
\end{equation}
where
\begin{footnotesize}\begin{equation*}
 A(z)=\begin{pmatrix}
        \lambda+\mu_1 & 0 & 0 & -\lambda p_1z & -\mu_2 & -\lambda p_1z & 0 & 0 & -\mu_3 & -\lambda p_1z & 0 & 0 \\[0.8em]
        \frac{-\lambda p_2}{p_2+p_3} & \lambda+\mu_1 & 0 & 0 & 0 & 0 & 0 & 0 & 0 & \frac{-\mu_3 q_1}{q_1+q_2} & 0 & 0 \\[0.8em]
        \frac{-\lambda p_3}{p_2+p_3} & 0 & \lambda+\mu_1 & 0 & 0 & \frac{-\mu_2q_1}{q_1+q_3} & 0 & 0 & 0 & 0 & 0 & 0 \\[0.8em]
        -\mu_1 & -\lambda z & -\lambda z & (\lambda+\mu_1)z & 0 & 0 & 0 & 0 & 0 & 0 & 0 & 0 \\[0.8em]
         0 & 0 & 0 & 0 & \lambda+\mu_2 & 0 & \frac{-\lambda zp_1}{p_1+p_3} & 0 & 0 & \frac{-\mu_3q_2}{q_1+q_2} & 0 & 0 \\[0.8em]
         0 & 0 & 0 & 0 & -\lambda & \lambda+\mu_2 & -\mu_2 & -\lambda z & 0 & 0 & 0 & 0 \\[0.8em]
         0 & -\mu_1 & 0 & -\lambda p_2z & 0 & -\lambda p_2z & (\lambda+\mu_2)z & 0 & 0 & -\lambda p_2z & 0 & -\mu_3z \\[0.8em]
         0 & 0 & 0 & \frac{-\mu_1q_2}{q_2+q_3} & 0 & 0 & \frac{-\lambda zp_3}{p_1+p_3} & (\lambda+\mu_2)z & 0 & 0 & 0 & 0 \\[0.8em]
         0 & 0 & 0 & 0 & 0 & \frac{-\mu_2q_3}{q_1+q_3} & 0 & 0 & \lambda+\mu_3 & 0 & \frac{-\lambda z p_1}{p_1+p_2} & 0 \\[0.8em]
         0 & 0 & 0 & 0 & 0 & 0 & 0 & 0 & -\lambda & \lambda+\mu_3 & -\mu_3 & -\lambda z \\[0.8em]
         0 & 0 & -\mu_1 & -\lambda p_3z & 0 & -\lambda p_3z & 0 & -\mu_2z & 0 & -\lambda p_3z & (\lambda+\mu_3)z & 0 \\[0.8em]
         0 & 0 & 0 & \frac{-\mu_1q_3}{q_2+q_3} & 0 & 0 & 0 & 0 & 0 & 0 & \frac{-\lambda zp_2}{p_1+p_2} & (\lambda+\mu_3)z
        \end{pmatrix},
\end{equation*}
\end{footnotesize}
\begin{align*}
\vec{G}(z)=\Big(G&_{(1;1,1)}(z),G_{(1;0,1)}(z),G_{(1;1,0)}(z),G_{(1;0,0)}(z),G_{(2;0,1)}(z),G_{(2;0,0)}(z),G_{(2;-1,0)}(z),\\
&G_{(2;-1,-1)}(z),G_{(3;1,0)}(z),G_{(3;0,0)}(z),G_{(3;0,-1)}(z),G_{(3;-1,-1)}(z)\Big)^T
\end{align*}
is a 12-dimensional column vector of the desired conditional PGFs, and
\begin{align*}
\vec{P}(z)=\Big(0&,\frac{-\mu_3q_1}{q_1+q_2}P_{0,(3;0,0)},\frac{-\mu_2q_1}{q_1+q_3}P_{0,(2;0,0)},\mu_1zP_{0,(1;0,0)},\frac{-\mu_3q_2}{q_1+q_2}P_{0,(3;0,0)},\mu_2P_{0,(2;0,0)},\\
&0,\frac{-\mu_1q_2}{q_2+q_3}P_{0,(1;0,0)},\frac{-\mu_2q_3}{q_1+q_3}P_{0,(2;0,0)},\mu_3P_{0,(3;0,0)},0,\frac{-\mu_1q_3}{q_2+q_3}P_{0,(1;0,0)}\Big)^T
\end{align*}
is a 12-dimensional vector containing the three unknowns, so-called 'boundary probabilities', $P_{0,(1;0,0)}$, $P_{0,(2;0,0)}$ and $P_{0,(3;0,0)}$.
Using Cramer's rule we have $G_d(z)=\frac{\mid A_d(z)\mid}{\mid A(z)\mid}$, $d\in\mathfrak{D}$,
where $|A|$ is the determinant of a matrix $A$, and $A_d(z)$ is the matrix obtained
from $A(z)$ by replacing the corresponding column of the latter matrix with $\vec{P}(z)$. Each of the PGFs $G_d(z)$, $d\in\mathfrak{D}$, is a function of the three unknown  probabilities, $P_{0,(1;0,0)}$, $P_{0,(2;0,0)}$ and $P_{0,(3;0,0)}$ appearing in $\vec{P}(z)$. In order to calculate those probabilities, we use first the normalization condition (\ref{SL3 eq sum=1}),
\begin{equation}\label{SL3 eq sum=1}
\sum_ {d\in\mathfrak{D}}G_d(1)=\sum_ {d\in\mathfrak{D}}\lim\limits_{z\rightarrow1}\frac{\mid A_d(z)\mid}{\mid A(z)\mid}=1,
\end{equation}
resulting in one equation in the above boundary probabilities, while two other equations between the three boundary probabilities
are derived as follows. Since $G_d(z)$ is defined for all $|z|\leq1$, each root of $|A(z)|$ is a root of $|A_d(z)|$. As a result, we utilize the two roots of $|A(z)|$ lying in $(-1,1)$, and solve $|A_d(z)|=0$, for arbitrary two different values of $d$. Once the boundary probabilities are calculated, the 12 PGFs can be obtained.

\subsection{Performance measures}\label{SL3 sec queue length}
In this section we calculate the mean and variance of $L_i$, $i=1,2,3$, as well as the covariances and correlation coefficients between each pair of queue lengths. We also derive the effective arrival rate to each queue, the probability that the server is idle, and the mean waiting time in each queue.

First
\begin{align*}
&\e[L_1]=\sum\limits_{d\in\mathfrak{D}}{G'_d(1)},\\
&\e[D_2]=\e[L_1-L_2]=G_{(1;1,1)}(1)+G_{(1;1,0)}(1)+G_{(3;1,0)}(1)-G_{(2;-1,0)}(1)-G_{(2;-1,-1)}(1)-G_{(3;-1,-1)}(1),\\
&\e[D_3]=\e[L_1-L_3]=G_{(1;1,1)}(1)+G_{(1;0,1)}(1)+G_{(2;0,1)}(1)-G_{(2;-1,-1)}(1)-G_{(3;0,-1)}(1)-G_{(3;-1,-1)}(1),\\
&\e[L_i]=\e[L_1]-\e[D_i], \ \ i=2,3.
\end{align*}
The above follows since the difference $L_1-L_2$ or $L_1-L_3$ can assume only the values -1, 0, or 1.
Next, for $j=2,3$,
\begin{align*}
Cov(L_1,L_j)&=\e[L_1L_j]-\e[L_1]\e[L_j]=\e[L_1(L_1-D_j)]-\e[L_1]\e[L_j]\\
&=\e[L_1^2]-\e[L_1D_j]-\e[L_1]\e[L_j],
\end{align*}
where
\begin{align*}
&\e[L_1^2]=\sum\limits_{d\in\mathfrak{D}}{G''_d(1)}+\e[L_1],\\
&\e[L_1D_2]=G'_{(1;1,1)}(1)+G'_{(1;1,0)}(1)+G'_{(3;1,0)}(1)-G'_{(2;-1,0)}(1)-G'_{(2;-1,-1)}(1)-G'_{(3;-1,-1)}(1),\\
&\e[L_1D_3]=G'_{(1;1,1)}(1)+G'_{(1;0,1)}(1)+G'_{(2;0,1)}(1)-G'_{(2;-1,-1)}(1)-G'_{(3;0,-1)}(1)-G'_{(3;-1,-1)}(1).\\
\end{align*}
Also,
\begin{align*}
Cov(L_2,L_3)&=\e[L_2L_3]-\e[L_2]\e[L_3]=\e[(L_1-D_2)(L_1-D_3)]-\e[L_2]\e[L_3]\\
&=\e[L_1^2]-\e[L_1D_2]-\e[L_1D_3]+\e[D_2D_3]-\e[L_2]\e[L_3],
\end{align*}
where
\[
\e[D_2D_3]=G_{(1;1,1)}(1)+G_{(2;-1,-1)}(1)+G_{(3;-1,-1)}(1).
\]

Furthermore, the variance of $L_i$, for $i=1,2,3$, is given by
\begin{align*}
&Var(L_1)=\e[L_1^2]-\left(\e[L_1]\right)^2,\\
&Var(L_j)=\e[L_j^2]-\left(\e[L_j]\right)^2=\e[(L_1-D_j)^2]-\left(\e[L_j]\right)^2\\
& \ \ \ \ \ \ \ \ \ \ \ =\e[L_1^2]-2\e[L_1D_j]+\e[D_j^2]-\left(\e[L_j]\right)^2,
\end{align*}
where
\begin{align*}
&\e[D_2^2]=G_{(1;1,1)}(1)+G_{(1;1,0)}(1)+G_{(2;-1,0)}(1)-G_{(2;-1,-1)}(1)-G_{(3;1,0)}(1)-G_{(3;-1,-1)}(1),\\
&\e[D_3^2]=G_{(1;1,1)}(1)+G_{(1;0,1)}(1)+G_{(2;0,1)}(1)-G_{(2;-1,-1)}(1)-G_{(3;0,-1)}(1)-G_{(3;-1,-1)}(1).
\end{align*}

From all the above, the correlation coefficient between $L_i$ and $L_j$, denoted by $Cor(L_i,L_j)$, can be explicitly calculated, using $Cor(L_i,L_j)=\frac{Cov(L_i,L_j)}{\sqrt{Var(L_i)Var(L_j)}}$.

Let $\lambda^i_{eff}$ denote the effective arrival rate to $Q_i$, $i=1,2,3$. Then,
\begin{align*}
\lambda^1_{eff}&=\lambda \Big(p_1(G_{(1;0,0)}(1)+G_{(2;0,0)}(1)+G_{(3;0,0)}(1))+G_{(2;-1,-1)}(1)+G_{(3;-1,-1)}(1)+\\
& \ \ \ \ \ \ \ \  \frac{p_1}{p_1+p_2}G_{(3;0,-1)}(1)+\frac{p_1}{p_1+p_3}G_{(2;-1,0)}(1)\Big),\\
\lambda^2_{eff}&=\lambda \Big(p_2(G_{(1;0,0)}(1)+G_{(2;0,0)}(1)+G_{(3;0,0)}(1))+G_{(1;1,0)}(1)+G_{(3;1,0)}(1)+\\
& \ \ \ \ \ \ \ \  \frac{p_2}{p_1+p_2}G_{(3;0,-1)}(1)+\frac{p_2}{p_2+p_3}G_{(1;1,1)}(1)\Big),\\
\lambda^3_{eff}&=\lambda \Big(p_3(G_{(1;0,0)}(1)+G_{(2;0,0)}(1)+G_{(3;0,0)}(1))+G_{(1;0,1)}(1)+G_{(2;0,1)}(1)+\\
& \ \ \ \ \ \ \ \  \frac{p_3}{p_1+p_3}G_{(2;-1,0)}(1)+\frac{p_3}{p_2+p_3}G_{(1;1,1)}(1)\Big).
\end{align*}
Clearly, $\sum_{i=1}^3{\lambda^i_{eff}}=\lambda$. Defining $\rho_{eff}^i=\frac{\lambda^i_{eff}}{\mu_i}$, we have that
\[
\p(\text{Server is idle})=P_{0,(1;0,0)}+P_{0,(2;0,0)}+P_{0,(3;0,0)}=1-\sum_{i=1}^3{\rho_{eff}^i}.
\]

Furthermore, define $W_i$ as the sojourn time of a customer in $Q_i$, for $i=1,2,3$. Then, by Little's Law,
\[\e[W_i]=\frac{\e[L_i]}{\lambda^i_{eff}}\ .\]
Last, denote by $\gamma_i$ the probability that the sever resides in $Q_i$, for $i=1,2,3$. Then,
\begin{align*}
&\gamma_1=G_{(1;1,1)}(1)+G_{(1;0,1)}(1)+G_{(1;1,0)}(1)+G_{(1;0,0)}(1),\\
&\gamma_2=G_{(2;0,1)}(1)+G_{(2;0,0)}(1)+G_{(2;-1,0)}(1)+G_{(2;-1,-1)}(1),\\
&\gamma_3=G_{(3;1,0)}(1)+G_{(3;0,0)}(1)+G_{(3;0,-1)}(1)+G_{(3;-1,-1)}(1).
\end{align*}

\section{Matrix Geometric Analysis}\label{SL3 sec: MG}
\subsection{Definitions and notations}\label{SL3 sec MG notations}
In this section we present an alternative approach to analyze the 3-queue JSQ-SLQ model by defining a Quasi Birth and Death (QBD) process, possessing 12 phases and infinite number of levels. With this analysis we derive the stability condition of the system. Phases are represented by the states $d$, $d\in\mathfrak{D}$, while levels are presented by the different values of $L_1\geq0$. For $n\geq 0$ define $\mathcal{S}_n$ to be the set of possible states for level $n$, where
\begin{align*}
\mathcal{S}_0=&\{(0,(1;0,0)),(0,(2;0,0)),(0,(2;-1,0)),(0,(2;-1,-1)),(0,(3;0,0)),(0,(3;0,-1)),(0,(3;-1,-1))\}, \\
\mathcal{S}_n=&\{(n,d)\mid d\in\mathfrak{D}\}, \ \   n\geq 1,
\end{align*}
and arrange the system's states in the order
\begin{align*}
\mathcal{S}=\Big\{\mathcal{S}_0; \mathcal{S}_1;\mathcal{S}_2;\ldots;\mathcal{S}_n\ldots\Big\}.
\end{align*}
The infinitesimal generator matrix of the QBD process, denoted by $Q$, is given by
\[
Q=\begin{pmatrix}
        B_1  & B_0 & 0 & \cdots & \cdots & \cdots \\[0.8em]
         B_2 & A_1 & A_0 & 0 & \cdots & \cdots \\[0.8em]
          0 &  A_2 & A_1 & A_0 & 0 & \cdots  \\[0.8em]
           \vdots & \ddots & \ddots & \ddots & \ddots & \ddots
        \end{pmatrix},
 \]
where,
\[
B_0=\begin{pmatrix}
        \lambda p_1  & 0 & 0 & 0 & 0  & 0 & 0 & 0 & 0  & 0 & 0 & 0 \\[0.3em]
        \lambda p_1  & 0 & 0 & 0 & 0  & 0 & 0 & 0 & 0  & 0 & 0 & 0 \\[0.3em]
        0  & 0 & 0 & 0 & \frac{\lambda p_1}{p_1+p_3}  & 0 & 0 & 0 & 0  & 0 & 0 & 0 \\[0.3em]
        0  & 0 & 0 & 0 & 0  & \lambda p_1 & 0 & 0 & 0  & 0 & 0 & 0 \\[0.3em]
        \lambda p_1  & 0 & 0 & 0 & 0  & 0 & 0 & 0 & 0  & 0 & 0 & 0 \\[0.3em]
        0  & 0 & 0 & 0 & 0  & 0 & 0 & 0 & \frac{\lambda p_1}{p_1+p_2}  & 0 & 0 & 0 \\[0.3em]
        0  & 0 & 0 & 0 & 0  & 0 & 0 & 0 & 0  & \lambda & 0 & 0
        \end{pmatrix}_{7\times12},
\]

\[B_1=\begin{pmatrix}
       -\lambda & 0 & \lambda p_2 & 0  & 0 & \lambda p_3 & 0 \\[0.3em]
       0 & -\lambda & \lambda p_2 & 0  & 0 & \lambda p_3 & 0 \\[0.3em]
       0 & \mu_2 & -(\lambda+\mu_2) & \frac{\lambda p_3}{p_1+p_3}  & 0 & 0 & 0 \\[0.3em]
       0 & 0 & 0 & -(\lambda+\mu_2)  & 0 & \mu_2 & 0 \\[0.3em]
       0 & 0 & \lambda p_2 & 0  & -\lambda & \lambda p_3 & 0 \\[0.3em]
       0 & 0 & 0 & 0 & \mu_3 & -(\lambda+\mu_3) & \frac{\lambda p_2}{p_1+p_2} \\[0.3em]
       0 & 0 & \mu_3 & 0  & 0 & 0 & -(\lambda+\mu_3)
        \end{pmatrix}_{7\times7},
\]

\[
B_2=\begin{pmatrix}
       \mu_1 & 0 & 0 & 0 & 0 & 0 & 0 \\[0.3em]
       0 & 0 & \mu_1 & 0 & 0 & 0 & 0 \\[0.3em]
       0 & 0 & 0 & 0 & 0 & \mu_1 & 0 \\[0.3em]
       0 & 0 & 0 & \frac{\mu_1q_2}{q_2+q_3} & 0 & 0 & \frac{\mu_1q_3}{q_2+q_3} \\[0.3em]
       0 & 0 & 0 & 0 & 0 & 0 & 0 \\[0.3em]
       0 & 0 & 0 & 0 & 0 & 0 & 0 \\[0.3em]
       0 & 0 & 0 & 0 & 0 & 0 & 0 \\[0.3em]
       0 & 0 & 0 & 0 & 0 & 0 & 0 \\[0.3em]
       0 & 0 & 0 & 0 & 0 & 0 & 0 \\[0.3em]
       0 & 0 & 0 & 0 & 0 & 0 & 0 \\[0.3em]
       0 & 0 & 0 & 0 & 0 & 0 & 0 \\[0.3em]
       0 & 0 & 0 & 0 & 0 & 0 & 0
        \end{pmatrix}_{12\times7}.
 \]
and
\[
A_0=\begin{pmatrix}
        0  & 0 & 0 & 0 & 0 & 0 & 0 & 0 & 0 & 0 & 0 & 0 \\[0.3em]
        0 & 0 & 0 & 0 & 0 & 0 & 0 & 0 & 0 & 0 & 0 & 0 \\[0.3em]
        0  & 0 & 0 & 0 & 0 & 0 & 0 & 0 & 0 & 0 & 0 & 0 \\[0.3em]
        \lambda p_1  & 0 & 0 & 0 & 0 & 0 & 0 & 0 & 0 & 0 & 0 & 0 \\[0.3em]
        0  & 0 & 0 & 0 & 0 & 0 & 0 & 0 & 0 & 0 & 0 & 0 \\[0.3em]
        \lambda p_1  & 0 & 0 & 0 & 0 & 0 & 0 & 0 & 0 & 0 & 0 & 0 \\[0.3em]
        0  & 0 & 0 & 0 & \frac{\lambda p_1}{p_1+p_3} & 0 & 0 & 0 & 0 & 0 & 0 & 0 \\[0.3em]
        0  & 0 & 0 & 0 & 0 & \lambda & 0 & 0 & 0 & 0 & 0 & 0 \\[0.3em]
        0  & 0 & 0 & 0 & 0 & 0 & 0 & 0 & 0 & 0 & 0 & 0 \\[0.3em]
         \lambda p_1  & 0 & 0 & 0 & 0 & 0 & 0 & 0 & 0 & 0 & 0 & 0 \\[0.3em]
        0  & 0 & 0 & 0 & 0 & 0 & 0 & 0 & \frac{\lambda p_1}{p_1+p_2} & 0 & 0 & 0 \\[0.3em]
        0  & 0 & 0 & 0 & 0 & 0 & 0 & 0 & 0 & \lambda & 0 & 0
        \end{pmatrix}_{12\times12},
\]

\begin{tiny}
\[
A_1=\begin{pmatrix}
        -(\lambda+\mu_1) & \frac{\lambda p_2}{p_2+p_3} & \frac{\lambda p_3}{p_2+p_3} & 0 & 0 & 0 & 0 & 0 & 0 & 0 & 0 & 0 \\[0.3em]
        0 & -(\lambda+\mu_1) & 0 & \lambda & 0 & 0 & 0 & 0 & 0 & 0 & 0 & 0 \\[0.3em]
        0 & 0 & -(\lambda+\mu_1) & \lambda & 0 & 0 & 0 & 0 & 0 & 0 & 0 & 0 \\[0.3em]
        0 & 0 & 0 & -(\lambda+\mu_1) & 0 & 0 & \lambda p_2 & 0 & 0 & 0 & \lambda p_3 & 0 \\[0.3em]
        \mu_2 & 0 & 0 & 0 & -(\lambda+\mu_2) & \lambda & 0 & 0 & 0 & 0 & 0 & 0 \\[0.3em]
        0 & 0 & \frac{\mu_2 q_1}{q_1+q_3} & 0 & 0 & -(\lambda+\mu_2) & \lambda p_2 & 0 & \frac{\mu_2 q_3}{q_1+q_3} & 0 & \lambda p_3 & 0 \\[0.3em]
        0 & 0 & 0 & 0 & 0 & \mu_2 & -(\lambda+\mu_2) & \frac{\lambda p_3}{p_1+p_3} & 0 & 0 & 0 & 0 \\[0.3em]
        0 & 0 & 0 & 0 & 0 & 0 & 0 & -(\lambda+\mu_2) & 0 & 0 & \mu_2 & 0 \\[0.3em]
        \mu_3 & 0 & 0 & 0 & 0 & 0 & 0 & 0 & -(\lambda+\mu_3) & \lambda & 0 & 0 \\[0.3em]
        0 & \frac{\mu_3 q_1}{q_1+q_2} & 0 & 0 & \frac{\mu_3 q_2}{q_1+q_2} & 0 & \lambda p_2 & 0 & 0 & -(\lambda+\mu_3) & \lambda p_3 & 0 \\[0.3em]
        0 & 0 & 0 & 0 & 0 & 0 & 0 & 0 & 0 & \mu_3 & -(\lambda+\mu_3) & \frac{\lambda p_2}{p_1+p_2} \\[0.3em]
        0 & 0 & 0 & 0 & 0 & 0 & \mu_3 & 0 & 0 & 0 & 0 & -(\lambda+\mu_3)
        \end{pmatrix}_{12\times12},
 \]
 \end{tiny}

\[
A_2=\begin{pmatrix}
      0  & 0 & 0 & \mu_1 & 0 & 0 & 0 & 0 & 0 & 0 & 0 & 0 \\[0.3em]
      0  & 0 & 0 & 0 & 0 & 0 & \mu_1 & 0 & 0 & 0 & 0 & 0 \\[0.3em]
      0  & 0 & 0 & 0 & 0 & 0 & 0 & 0 & 0 & 0 & \mu_1 & 0 \\[0.3em]
      0  & 0 & 0 & 0 & 0 & 0 & 0 & \frac{\mu_1 q_2}{q_2+q_3} & 0 & 0 & 0 & \frac{\mu_1 q_3}{q_2+q_3} \\[0.3em]
      0  & 0 & 0 & 0 & 0 & 0 & 0 & 0 & 0 & 0 & 0 & 0 \\[0.3em]
      0  & 0 & 0 & 0 & 0 & 0 & 0 & 0 & 0 & 0 & 0 & 0 \\[0.3em]
      0  & 0 & 0 & 0 & 0 & 0 & 0 & 0 & 0 & 0 & 0 & 0 \\[0.3em]
      0  & 0 & 0 & 0 & 0 & 0 & 0 & 0 & 0 & 0 & 0 & 0 \\[0.3em]
      0  & 0 & 0 & 0 & 0 & 0 & 0 & 0 & 0 & 0 & 0 & 0 \\[0.3em]
      0  & 0 & 0 & 0 & 0 & 0 & 0 & 0 & 0 & 0 & 0 & 0 \\[0.3em]
      0  & 0 & 0 & 0 & 0 & 0 & 0 & 0 & 0 & 0 & 0 & 0 \\[0.3em]
      0  & 0 & 0 & 0 & 0 & 0 & 0 & 0 & 0 & 0 & 0 & 0
        \end{pmatrix}_{12\times12}.
\]

\subsection{Stability condition}\label{SL3 sec MG stability}
Define the matrix $A=A_0+A_1+A_2$. We get
\begin{tiny}
\[
A=\begin{pmatrix}
        -(\lambda+\mu_1) & \frac{\lambda p_2}{p_2+p_3} & \frac{\lambda p_3}{p_2+p_3} & \mu_1 & 0 & 0 & 0 & 0 & 0 & 0 & 0 & 0 \\[0.3em]
        0 & -(\lambda+\mu_1) & 0 & \lambda & 0 & 0 & \mu_1 & 0 & 0 & 0 & 0 & 0 \\[0.3em]
        0 & 0 & -(\lambda+\mu_1) & \lambda & 0 & 0 & 0 & 0 & 0 & 0 & \mu_1 & 0 \\[0.3em]
        \lambda p_1 & 0 & 0 & -(\lambda+\mu_1) & 0 & 0 & \lambda p_2 & \frac{\mu_1 q_2}{q_2+q_3} & 0 & 0 & \lambda p_3 & \frac{\mu_1 q_3}{q_2+q_3} \\[0.3em]
        \mu_2 & 0 & 0 & 0 & -(\lambda+\mu_2) & \lambda & 0 & 0 & 0 & 0 & 0 & 0 \\[0.3em]
        \lambda p_1 & 0 & \frac{\mu_2 q_1}{q_1+q_3} & 0 & 0 & -(\lambda+\mu_2) & \lambda p_2 & 0 & \frac{\mu_2 q_3}{q_1+q_3} & 0 & \lambda p_3 & 0 \\[0.3em]
        0 & 0 & 0 & 0 & \frac{\lambda p_1}{p_1+p_3} & \mu_2 & -(\lambda+\mu_2) & \frac{\lambda p_3}{p_1+p_3} & 0 & 0 & 0 & 0 \\[0.3em]
        0 & 0 & 0 & 0 & 0 & \lambda & 0 & -(\lambda+\mu_2) & 0 & 0 & \mu_2 & 0 \\[0.3em]
        \mu_3 & 0 & 0 & 0 & 0 & 0 & 0 & 0 & -(\lambda+\mu_3) & \lambda & 0 & 0 \\[0.3em]
        \lambda p_1 & \frac{\mu_3 q_1}{q_1+q_2} & 0 & 0 & \frac{\mu_3 q_2}{q_1+q_2} & 0 & \lambda p_2 & 0 & 0 & -(\lambda+\mu_3) & \lambda p_3 & 0 \\[0.3em]
        0 & 0 & 0 & 0 & 0 & 0 & 0 & 0 & \frac{\lambda p_1}{p_1+p_2} & \mu_3 & -(\lambda+\mu_3) & \frac{\lambda p_2}{p_1+p_2} \\[0.3em]
        0 & 0 & 0 & 0 & 0 & 0 & \mu_3 & 0 & 0 & \lambda & 0 & -(\lambda+\mu_3)
        \end{pmatrix}_{12\times12},
\]
\end{tiny}
The matrix $A$ is the infinitesimal generator matrix of the process describing the evolution of $D$, given that $L_1\geq 1$. Let $\vec{\pi}=\left(\pi_1,\pi_2,...,\pi_{12}\right)$ be the stationary vector of the matrix $A$, i.e. $\vec{\pi}A=\vec{0}$ and $\vec{\pi}\cdot\vec{e}=1$ (where $\vec{e}$  is a 12-dimensional column vector with all its entries equal to 1). From \cite{neuts1981matrix}, we have that the stability condition is
\[
\vec{\pi}A_0\vec{e}<\vec{\pi}A_2\vec{e}.
\]
\cite{hanukov2020explicit} showed that if the 3 matrices $A_0$, $A_1$ and $A_2$ are all upper triangular, or all lower triangular, then there exists
a simple and direct stability condition. Unfortunately, this is not the case here, so we have to use the above condition, which translates into
\[
\lambda\left(\pi_8+\pi_{12}+p_1(\pi_4+\pi_6+\pi_{10})+\frac{p_1}{p_1+p_2}\pi_{11}+\frac{p_1}{p_1+p_3}\pi_7\right)<\mu_1(\pi_1+\pi_2+\pi_3+\pi_4).
\]
Explicit numerical expressions for the elements of $\vec{\pi}$ can be calculated by any mathematical software. We omit these expressions. However, as an illustration, we assume values for the parameters $\mu_i$, $p_i$ and $q_i$ for $i=1,2,3$, and define
\[
f(\lambda)=\lambda\left(\pi_8+\pi_{12}+p_1(\pi_4+\pi_6+\pi_{10})+\frac{p_1}{p_1+p_2}\pi_{11}+\frac{p_1}{p_1+p_3}\pi_7\right)-\mu_1(\pi_1+\pi_2+\pi_3+\pi_4).
\]
That is, the system is stable iff $f(\lambda)<0$. For example, consider the symmetric case where $\mu_i=5$, $p_i=\frac{1}{3}$ and $q_i=1$ for $i=1,2,3$. We then get
\[
\vec{\pi}=\left(\frac{1}{9},\frac{1}{18},\frac{1}{18},\frac{1}{9},\frac{1}{18},\frac{1}{9},\frac{1}{9},\frac{1}{18},\frac{1}{18},\frac{1}{9},\frac{1}{9},\frac{1}{18}\right),
\]
while the stability condition $f(\lambda)<0$ translates into $\frac{\lambda-5}{3}<0$, so the system is stable, as expected, iff $\lambda<5$. Further examples for the stability condition in terms of the arrival rate $\lambda$ for various sets of parameters are given in Table \ref{SL3 stability examples}.
\begin{table}[h]
 \caption{Stability condition in terms of $\lambda$.}
  \label{SL3 stability examples}
    \centering\footnotesize{
\begin{tabular}{|l|l|}
  \hline
  Parameters & Stability condition \\ \hline
  $\mu_1=4$, $\mu_2=5$, $\mu_3=6$, $p_i=\frac{1}{3}$ and $q_i=1$ for $i=1,2,3$  & $\lambda<4.90215$ \\
  $\mu_1=4$, $\mu_2=5$, $\mu_3=6$, $p_1=0.2$, $p_2=0.3$, $p_3=0.5$ and $q_i=1$ for $i=1,2,3$ & $\lambda<5.0149$ \\
  $\mu_1=4$, $\mu_2=5$, $\mu_3=6$, $p_1=0.05$, $p_2=0.05$, $p_3=0.9$ and $q_i=1$ for $i=1,2,3$ & $\lambda<5.17335$ \\
  $\mu_1=4$, $\mu_2=5$, $\mu_3=6$, $p_1=0.9$, $p_2=0.05$, $p_3=0.05$ and $q_i=1$ for $i=1,2,3$ & $\lambda<4.63355$ \\ \hline
  \end{tabular}}
\end{table}

\subsection{Calculation of the system's states steady state distribution}\label{SL3 sec MG equilibrium}
For $n\geq 0$ define the steady-state probability vector $\vec{P}_n$, as follows:
\begin{equation*}
\vec{P}_n=\left\{ \begin{array}{ll}
\left(P_{0,(1;0,0)},P_{0,(2;0,0)},P_{0,(2;-1,0)},P_{0,(2;-1,-1)},P_{0,(3;0,0)},P_{0,(3;0,-1)},P_{0,(3;-1,-1)}\right),  \ n=0,  \\
\left(P_{n,(1;0,0)},P_{n,(1;0,1)},\ldots, P_{n,(3;-1,-1)}\right), \ \ n\geq 1.  \\
\end{array} \right.
\end{equation*}

From \cite{neuts1981matrix},
\[\vec{P}_n=\vec{P}_1R^{n-1}, \ \ n\geq 1,\]
where $R$ is the minimal non-negative solution of the matrix quadratic equation
\begin{equation}\label{SL eq: MG equation}
A_0+RA_1+R^2A_2=0.
\end{equation}
The matrix $R$ can be calculated by using well-known algorithms, see e.g. \cite{neuts1981matrix}, \cite{Latouche}, \cite{Artalejo}, and \cite{harchol2013performance}.
The vectors $\vec{P}_0$ and $\vec{P}_1$  are the solutions of the following linear system of equations:
\begin{align}
&\vec{P}_0B_1+\vec{P}_1B_2=\vec{0}, \nonumber \\
&\vec{P}_{0}B_0+\vec{P}_1(A_1+RA_2)=\vec{0} ,  \nonumber \\
&\vec{P}_0\vec{e}_0+\vec{P}_{1}[\textbf{I}-R]^{-1}\vec{e}=1,\label{SL3non eq normalization}
\end{align}
where $\vec{e}_0$ is a $7$-dimensional vector of 1's. Equation (\ref{SL3non eq normalization}) is the normalization equation.


\section{Numerical Results and the Gini Index}\label{SL3 sec results}
In this section we present numerical results for the performance measures derived in Section \ref{SL3 sec queue length}. We calculate $\e[L_i]$, $\e[W_i]$, $i=1,2,3$, for a wide range of parameters, as well as the correlation coefficients $Cor(L_i,L_j)$ between each pair of queue lengths. We also calculate $\lambda^i_{eff}$, i.e. the effective arrival rate to queue $i$, as well as $\gamma_i$, which is the probability that the server resides in $Q_i$, $i=1,2,3$. Tables \ref{SL3 mu2=3 mu3=5 p1=p2=p3=1/3 and q1=q2=q3=1}-\ref{SL3 gamma mu2=5, mu3=5 p1=1, p2=p3=0 and q1=q2=q3=1} examine the influence of $\mu_1$ on the various performance measures. Furthermore, we calculate the Gini coefficient, (see e.g. \cite[pp.~231--233]{dodge2008concise}), which is a well known measure for inequality. The Gini coefficient with respect to the differences between the mean queue sizes in the JSQ-SLQ operating model is calculated as follows:
\[
GI=\frac{\sum\limits_{i=1}^{3}{\sum\limits_{i=1}^{3}{\big|\e[L_i]-\e[L_j]\big|}}}{2\cdot3\cdot\sum\limits_{i=1}^{3}{\e[L_i]}}.
\]
It turns out, as it is shown in the tables below, that the Gini coefficient is close to zero in almost all sets of parameters, even when the system is close to saturation (resulting in large mean queue sizes). This outcome demonstrates the effectiveness of the combined JSQ-SLQ policy in balancing the queue sizes.

In Tables \ref{SL3 mu2=3 mu3=5 p1=p2=p3=1/3 and q1=q2=q3=1}-\ref{SL3 gamma mu2=5, mu3=5, p1=p2=p3=1/3 and q1=q2=q3=1} we set $\lambda=4$, $\mu_3=5$, $p_1=p_2=p_3=\frac{1}{3}$ and $q_1=q_2=q_3=1$, while $\mu_2$ varies between the tables (assuming the values 3,4,5), and $\mu_1$ varies within each table (where its minimum is close to the value ensuring stability). In Tables \ref{SL3 mu2=3, mu3=5 p1=1, p2=p3=0 and q1=q2=q3=1}-\ref{SL3 gamma mu2=5, mu3=5 p1=1, p2=p3=0 and q1=q2=q3=1} we set $\lambda=4$, $\mu_3=5$, $p_1=0.9999$, $p_2=p_3=0.00005$, $q_1=q_2=q_3=1$, while $\mu_2$ again varies between the tables assuming the values 3, 4 and 5. Tables \ref{SL3 q1 mu1=4.5 mu2=3, mu3=5 p1=p2=p3=1/3 and q2=q3=1}-\ref{SL3 gamma q1 mu1=4.2 mu2=3, mu3=5 p1=1, p2=p3=0 and q2=q3=1} examine the impact of $q_1$, which is the weight assigned to $Q_1$ from the server's point of view.

\begin{table}[h]
 \caption{Numerical results for $\e[L_i]$, $\e[W_i]$ and $Cor(L_i,L_j)$ when $\lambda=4$, $\mu_2=3$, $\mu_3=5$, $p_1=p_2=p_3=\frac{1}{3}$, $q_1=q_2=q_3=1$}
  \label{SL3 mu2=3 mu3=5 p1=p2=p3=1/3 and q1=q2=q3=1}
    \centering\footnotesize{
\begin{tabular}{|l||c|c|c|c|c|c|c|c|c|c|c|}
  \hline
  $\mu_1$  & $\e[L_1]$ & $\e[L_2]$ &$\e[L_3]$& $\e[W_1]$ & $\e[W_2]$& $\e[W_3]$ & $Cor(L_1,L_2)$ & $Cor(L_1,L_3)$ & $Cor(L_2,L_3)$ & $GI$ \\ \hline
  $\mu_1=4.5$ & 48.64 & 48.69 & 48.62 & 35.46 & 39.90 & 34.53 &  0.9999 & 0.9999 &0.9999 & 0.0003\\
   $\mu_1=5$ & 8.31 & 8.37 & 8.31 & 5.96 & 6.91 & 5.96 &  0.9970 & 0.9970 &0.9970 & 0.0016\\
   $\mu_1=10$ & 1.33 & 1.46 & 1.40 & 0.88 & 1.25 & 1.05 &  0.9266 & 0.9260 &0.9219 & 0.0207\\
  $\mu_1=30$ & 0.67 & 0.85 & 0.79 & 0.42 & 0.76 & 0.62 &  0.8438 & 0.8390 & 0.8119 & 0.0519\\  \hline
\end{tabular}}
\end{table}

\begin{table}[h]
 \caption{Numerical results for $\lambda^i_{eff}$ and $\gamma_i$ when $\lambda=4$, $\mu_2=3$, $\mu_3=5$, $p_1=p_2=p_3=\frac{1}{3}$, $q_1=q_2=q_3=1$}
  \label{SL3 gamma mu2=3 mu3=5 p1=p2=p3=1/3 and q1=q2=q3=1}
    \centering\footnotesize{
\begin{tabular}{|l||c|c|c|c|c|c|}
  \hline
  $\mu_1$  & $\lambda^1_{eff}$ & $\lambda^2_{eff}$ & $\lambda^3_{eff}$ & $\gamma_1$ & $\gamma_2$ & $\gamma_3$ \\ \hline
  $\mu_1=4.5$ & 1.37 & 1.22 & 1.41 & 0.307 & 0.409 & 0.284  \\
   $\mu_1=5$ & 1.395 & 1.21 & 1.395 & 0.293 & 0.414 & 0.293  \\
   $\mu_1=10$ & 1.51 & 1.16 & 1.33 & 0.232 & 0.438 & 0.33  \\
  $\mu_1=30$ & 1.60 & 1.12 & 1.28 & 0.204 & 0.446 & 0.350  \\  \hline
\end{tabular}}
\end{table}

\begin{table}[h]
 \caption{Numerical results for $\e[L_i]$, $\e[W_i]$ and $Cor(L_i,L_j)$ when $\lambda=4$, $\mu_2=4$, $\mu_3=5$, $p_1=p_2=p_3=\frac{1}{3}$, $q_1=q_2=q_3=1$}
  \label{SL3 mu2=4, mu3=5,  p1=p2=p3=1/3 and q1=q2=q3=1}
    \centering\footnotesize{
\begin{tabular}{|l||c|c|c|c|c|c|c|c|c|c|}
  \hline
  $\mu_1$  & $\e[L_1]$ & $\e[L_2]$ &$\e[L_3]$& $\e[W_1]$ & $\e[W_2]$& $\e[W_3]$ & $Cor(L_1,L_2)$ & $Cor(L_1,L_3)$ & $Cor(L_2,L_3)$ & $GI$ \\ \hline
  $\mu_1=3.3$ & 88.46 & 88.44 & 88.41 & 70.28 & 66.44 & 62.69 &  0.9999 & 0.9999 &0.9999 & 0.0001\\
   $\mu_1=3.4$ & 22.85 & 22.83 & 22.80 & 18.04 & 17.20 & 16.22 &  0.9996 & 0.9996 &0.9996 & 0.0005\\
   $\mu_1=3.5$ & 13.26 & 13.25 & 13.22 & 10.41 & 10.01 & 9.43 &  0.9988 & 0.9988 &0.9988 & 0.0007\\
  $\mu_1=4$ & 4.49 & 4.49 & 4.47 & 3.44 & 3.44 & 3.23 &  0.9902 & 0.9902 & 0.9902 & 0.0010\\
  $\mu_1=4.5$ & 2.82 & 2.84 & 2.81 & 2.11 & 2.19 & 2.05 &  0.9768 & 0.9769 & 0.9769 & 0.0024\\
  $\mu_1=5$ & 2.11 & 2.14 & 2.11 & 1.56 & 1.67 & 1.56 &  0.9616 & 0.9618 & 0.9616 & 0.0031\\
  $\mu_1=10$ & 0.80 & 0.89 & 0.86 & 0.54 & 0.72 & 0.66 &  0.8434 & 0.8434 & 0.8341 & 0.0235\\
  $\mu_1=30$ & 0.45 & 0.59 & 0.56 & 0.29 & 0.50 & 0.45 &  0.7542 & 0.75154 & 0.7069 & 0.0583\\  \hline
\end{tabular}}
\end{table}

\begin{table}[h]
 \caption{Numerical results for $\lambda^i_{eff}$ and $\gamma_i$ when $\lambda=4$, $\mu_2=4$, $\mu_3=5$, $p_1=p_2=p_3=\frac{1}{3}$, $q_1=q_2=q_3=1$}
  \label{SL3 gamma mu2=4, mu3=5,  p1=p2=p3=1/3 and q1=q2=q3=1}
    \centering\footnotesize{
\begin{tabular}{|l||c|c|c|c|c|c|}
  \hline
  $\mu_1$  & $\lambda^1_{eff}$ & $\lambda^2_{eff}$ & $\lambda^3_{eff}$ & $\gamma_1$ & $\gamma_2$ & $\gamma_3$ \\ \hline
  $\mu_1=3.3$ & 1.26 & 1.33 & 1.41 & 0.383 & 0.334 & 0.283  \\
   $\mu_1=3.4$ & 1.27 & 1.33 & 1.40 & 0.377 & 0.337 & 0.286  \\
   $\mu_1=3.5$ & 1.27 & 1.32 & 1.41 & 0.372 & 0.339 & 0.289  \\
  $\mu_1=4$ & 1.31 & 1.31 & 1.38 & 0.349 & 0.349 & 0.302  \\
  $\mu_1=4.5$ & 1.34 & 1.29 & 1.37 & 0.332 & 0.357 & 0.311 \\
  $\mu_1=5$ & 1.36 & 1.28 & 1.36 & 0.318 & 0.364 & 0.318  \\
  $\mu_1=10$ & 1.47 & 1.23 & 1.30 & 0.262 & 0.389 & 0.349  \\
  $\mu_1=30$ & 1.56 & 1.19 & 1.25 & 0.237 & 0.399 & 0.364  \\  \hline
\end{tabular}}
\end{table}

\begin{table}[h]
 \caption{Numerical results for $\e[L_i]$, $\e[W_i]$ and $Cor(L_i,L_j)$ when $\lambda=4$, $\mu_2=5$, $\mu_3=5$, $p_1=p_2=p_3=\frac{1}{3}$, $q_1=q_2=q_3=1$}
  \label{SL3 mu2=5, mu3=5, p1=p2=p3=1/3 and q1=q2=q3=1}
    \centering\footnotesize{
\begin{tabular}{|l||c|c|c|c|c|c|c|c|c|c|}
  \hline
  $\mu_1$ & $\e[L_1]$ & $\e[L_2]$ &$\e[L_3]$& $\e[W_1]$ & $\e[W_2]$& $\e[W_3]$ & $Cor(L_1,L_2)$ & $Cor(L_1,L_3)$ & $Cor(L_2,L_3)$ & $GI$\\ \hline
  $\mu_1=2.75$ & 55.11 & 55.03 & 55.03 & 46.49 & 39.10 & 39.10 &  0.9999 & 0.9999 &0.9999 & 0.0003\\
   $\mu_1=3$ & 8.37 & 8.31 & 8.31 & 6.91 & 5.95 & 5.95 &  0.9970 & 0.9970 &0.9970 & 0.0016\\
    $\mu_1=3.4$ & 3.73 & 3.68 & 3.68 & 2.993 & 2.67 & 2.67 &  0.9859 & 0.9859 & 0.9859 & 0.0030\\
   $\mu_1=3.5$ & 3.29 & 3.25 & 3.25 & 2.63 & 2.37 & 2.37 &  0.9823 & 0.9823 &0.9824 & 0.0027\\
  $\mu_1=4$ & 2.14 & 2.11 & 2.11 & 1.67 & 1.56 & 1.56 &  0.9616 & 0.9616 & 0.9618 & 0.0031\\
  $\mu_1=4.5$ & 1.63 & 1.61 & 1.61 & 1.24 & 1.2 & 1.2 &  0.9389 & 0.9389 & 0.9391 & 0.0027\\
  $\mu_1=5$ & 1.33 & 1.33 & 1.33 & 1 & 1 & 1 &  0.9164 & 0.9164 & 0.9164 & 0\\
  $\mu_1=10$ & 0.60 & 0.66 & 0.66 & 0.42 & 0.52 & 0.52 &  0.7744 & 0.7744 & 0.7624 & 0.0208 \\
  $\mu_1=30$ & 0.35 & 0.46 & 0.46 & 0.22 & 0.38 & 0.38 &  0.6870 & 0.6870 & 0.6316 & 0.0577\\  \hline
\end{tabular}}
\end{table}

\begin{table}[h]
 \caption{Numerical results for $\lambda^i_{eff}$ and $\gamma_i$ when $\lambda=4$, $\mu_2=5$, $\mu_3=5$, $p_1=p_2=p_3=\frac{1}{3}$, $q_1=q_2=q_3=1$}
  \label{SL3 gamma mu2=5, mu3=5,  p1=p2=p3=1/3 and q1=q2=q3=1}
    \centering\footnotesize{
\begin{tabular}{|l||c|c|c|c|c|c|}
  \hline
  $\mu_1$  & $\lambda^1_{eff}$ & $\lambda^2_{eff}$ & $\lambda^3_{eff}$ & $\gamma_1$ & $\gamma_2$ & $\gamma_3$ \\ \hline
  $\mu_1=2.75$ & 1.18 & 1.41 & 1.41 & 0.432 & 0.284 & 0.284  \\
  $\mu_1=3$ & 1.21 & 1.39 & 1.39 & 0.414 & 0.293 & 0.293  \\
   $\mu_1=3.4$ & 1.24 & 1.38 & 1.38 & 0.390 & 0.305 & 0.305  \\
   $\mu_1=3.5$ & 1.26 & 1.37 & 1.37 & 0.386 & 0.307 & 0.307  \\
  $\mu_1=4$ & 1.28 & 1.36 & 1.36 & 0.364 & 0.318 & 0.318  \\
  $\mu_1=4.5$ & 1.31 & 1.345 & 1.345 & 0.347 & 0.3265 & 0.3265 \\
  $\mu_1=5$ & $\frac{4}{3}$ & $\frac{4}{3}$ & $\frac{4}{3}$ & $\frac{1}{3}$ & $\frac{1}{3}$ & $\frac{1}{3}$  \\
  $\mu_1=10$ & 1.44 & 1.28 & 1.28 & 0.280 & 0.360 & 0.360  \\
  $\mu_1=30$ & 1.53 & 1.235 & 1.235 & 0.257 & 0.3715 & 0.3715  \\  \hline
\end{tabular}}
\end{table}

The following insights arise from the tables:
\begin{enumerate}
\item In all tables, $\e[L_i]$ and $\e[W_i]$ ($i=1,2,3$) are all decreasing functions of $\mu_1$.
\item The correlation coefficient between $L_i$ and $L_j$ is always positive, since customers always join the shortest queue. Furthermore, as $\mu_1$ increases, the correlation coefficient between $L_i$ and $L_j$ decreases.  This phenomenon occurs since for small values of $\mu_1$, each service in $Q_1$ lasts for a long time and therefore, as a result of the JSQ policy, the number of customers in all queues increases simultaneously. On the other hand, for large values of $\mu_1$, the behavior of $L_j$ (for $j=2,3$) is less affected by $L_1$, since the server resides in $Q_1$ a short amount of time.
\item The Gini coefficient increases monotonically when $\mu_1$ increases, but remains very small (and close to zero is some cases). It is higher when $\mu_1$ is large. This follows since, even if the mean queue sizes are small, the relative differences between them are more significant.
\item When $p_1$ is very close to 1, the expected sojourn time in $Q_1$ is smaller than the expected sojourn times in $Q_2$ and $Q_3$. This occurs since when $p_1$ approaches 1, $\lambda^1_{eff}$ increases, which, due to the SLQ regime, increases the probability that the server resides in $Q_1$. Clearly, the presence of the server in $Q_1$ reduces the waiting time there.
\item When the value of $p_1$ approaches 1, the minimum value of $\mu_1$ that ensures the system's stability rises. For example, in Table \ref{SL3 mu2=5, mu3=5, p1=p2=p3=1/3 and q1=q2=q3=1}, the minimal value of $\mu_1$ for stability is 2.75, while in Table \ref{SL3 mu2=5, mu3=5 p1=1, p2=p3=0 and q1=q2=q3=1} it is 3.4.
\item In all tables, as $\mu_1$ increases, the probability that the server resides in $Q_1$ decreases.
\item Tables \ref{SL3 q1 mu1=4.5 mu2=3, mu3=5 p1=p2=p3=1/3 and q2=q3=1} and \ref{SL3 q1 mu1=4.5 mu2=3, mu3=5 p1=1, p2=p3=0 and q2=q3=1} show that The correlation coefficient between $L_i$ and $L_j$ remains unchanged as $q_1$ increases, while $\e[L_i]$ and $\e[W_i]$ decrease. Furthermore, Tables \ref{SL3 gamma q1 mu1=4.5 mu2=3, mu3=5 p1=p2=p3=1/3 and q2=q3=1} and \ref{SL3 gamma q1 mu1=4.2 mu2=3, mu3=5 p1=1, p2=p3=0 and q2=q3=1} show that as $q_1$ increases, the probability that the server resides in $Q_1$ slightly increases, while the probabilities that the server resides in $Q_2$ and $Q_3$  moderately decrease.
\item The numerical results strongly demonstrate that the aim of the combined policy JSQ+SLQ is very well achieved: in all tables, the mean queue sizes are almost equal, as well as the mean sojourn times.
\end{enumerate}

\begin{table}[h]
 \caption{Numerical results for $\e[L_i]$, $\e[W_i]$ and $Cor(L_i,L_j)$ when $\lambda=4$, $\mu_2=3$, $\mu_3=5$, $p_1=0.9999$, $p_2=p_3=0.00005$, $q_1=q_2=q_3=1$}
  \label{SL3 mu2=3, mu3=5 p1=1, p2=p3=0 and q1=q2=q3=1}
    \centering\footnotesize{
\begin{tabular}{|l||c|c|c|c|c|c|c|c|c|c|}
  \hline
  $\mu_1$  & $\e[L_1]$ & $\e[L_2]$ &$\e[L_3]$& $\e[W_1]$ & $\e[W_2]$& $\e[W_3]$ & $Cor(L_1,L_2)$ & $Cor(L_1,L_3)$ & $Cor(L_2,L_3)$ & $GI$ \\ \hline
  $\mu_1=4.2$ & 105.10 & 104.95 & 104.90 & 47.79 & 123.78 & 110.07 &  0.9999 & 0.9999 &0.9999 & 0.0004 \\
   $\mu_1=4.5$& 8.19 & 8.04 & 8.00 & 3.65 & 9.74 & 8.62 &  0.9967 & 0.9967 &0.9974 & 0.0052 \\
   $\mu_1=5$ & 3.36 & 3.22 & 3.17 & 1.45 & 4.06 & 3.56 &  0.9814 & 0.9819 &0.9858 & 0.0130\\
   $\mu_1=10$ & 0.63 & 0.52 & 0.48 & 0.23 & 0.90 & 0.75 &  0.7225 & 0.7250 & 0.7910 & 0.0614\\
  $\mu_1=30$ & 0.20 & 0.14 & 0.12 & 0.06 & 0.51 & 0.39 &  0.4698 & 0.4558 & 0.5629 & 0.1159\\  \hline
\end{tabular}}
\end{table}

\begin{table}[h]
 \caption{Numerical results for $\lambda^i_{eff}$ and $\gamma_i$ when $\lambda=4$, $\mu_2=3$, $\mu_3=5$, $p_1=0.9999$, $p_2=p_3=0.00005$, $q_1=q_2=q_3=1$}
  \label{SL3 gamma mu2=3, mu3=5 p1=1, p2=p3=0 and q1=q2=q3=1}
    \centering\footnotesize{
\begin{tabular}{|l||c|c|c|c|c|c|c|}
  \hline
  $\mu_1$  & $\lambda^1_{eff}$ & $\lambda^2_{eff}$ & $\lambda^3_{eff}$ & $\gamma_1$ & $\gamma_2$ & $\gamma_3$ \\ \hline
  $\mu_1=4.2$ & 2.20 & 0.85 & 0.95 & 0.526 & 0.283 & 0.191  \\
   $\mu_1=4.5$& 2.24 & 0.83 & 0.93 & 0.526 & 0.283 & 0.191  \\
   $\mu_1=5$ & 2.32 & 0.79 & 0.89 & 0.528 & 0.278 & 0.194    \\
   $\mu_1=10$ & 2.78 & 0.58 & 0.64 & 0.602 & 0.227& 0.171  \\
  $\mu_1=30$ & 3.41 & 0.28 & 0.31 & 0.790 & 0.117 & 0.093  \\  \hline
\end{tabular}}
\end{table}

\begin{table}[h]
 \caption{Numerical results for $\e[L_i]$, $\e[W_i]$ and $Cor(L_i,L_j)$ when $\lambda=4$, $\mu_2=4$, $\mu_3=5$, $p_1=0.9999$, $p_2=p_3=0.00005$, $q_1=q_2=q_3=1$}
  \label{SL3 mu2=4, mu3=5 p1=1, p2=p3=0 and q1=q2=q3=1}
    \centering\footnotesize{
\begin{tabular}{|l||c|c|c|c|c|c|c|c|c|c|}
  \hline
  $\mu_1$  & $\e[L_1]$ & $\e[L_2]$ &$\e[L_3]$& $\e[W_1]$ & $\e[W_2]$& $\e[W_3]$ & $Cor(L_1,L_2)$ & $Cor(L_1,L_3)$ & $Cor(L_2,L_3)$ & $GI$ \\ \hline
  $\mu_1=3.7$ & 59.77 & 59.57 & 59.55 & 28.32 & 64.51 & 61.63 &  0.9999 & 0.9999 &0.9999 & 0.0008\\
   $\mu_1=4$& 6.89 & 6.69 & 6.67 & 3.18 & 7.46 & 7.11 &  0.9954 & 0.9954 &0.9963 & 0.0072\\
   $\mu_1=4.5$ & 3.42 & 3.20 & 3.20 & 1.58 & 3.48 & 3.48 &  0.9824 & 0.9824 &0.9860 & 0.0149 \\
   $\mu_1=5$ & 2.91 & 2.72 & 2.70 & 1.29 & 3.17 & 3.02 &  0.9757 & 0.9760& 0.9811 & 0.0168 \\
  $\mu_1=10$ & 0.53 & 0.39 & 0.38 & 0.19 & 0.67 & 0.62 &  0.6530 & 0.6547 & 0.7349 & 0.0769 \\
  $\mu_1=30$ & 0.18 & 0.11 & 0.10 & 0.05 & 0.39 & 0.35 &  0.4146 & 0.4083 & 0.5083 & 0.1367 \\  \hline
\end{tabular}}
\end{table}

\begin{table}[h]
 \caption{Numerical results for $\lambda^i_{eff}$ and $\gamma_i$ when $\lambda=4$, $\mu_2=4$, $\mu_3=5$, $p_1=0.9999$, $p_2=p_3=0.00005$, $q_1=q_2=q_3=1$}
  \label{SL3 gamma mu2=4, mu3=5 p1=1, p2=p3=0 and q1=q2=q3=1}
    \centering\footnotesize{
\begin{tabular}{|l||c|c|c|c|c|c|c|}
  \hline
  $\mu_1$  & $\lambda^1_{eff}$ & $\lambda^2_{eff}$ & $\lambda^3_{eff}$ & $\gamma_1$ & $\gamma_2$ & $\gamma_3$ \\ \hline
  $\mu_1=3.7$ & 2.11 & 0.92 & 0.97 & 0.574 & 0.231 & 0.195 \\
   $\mu_1=4$& 2.17 & 0.90 & 0.93 & 0.571 & 0.232 & 0.197  \\
   $\mu_1=4.5$& 2.25 & 0.86 & 0.89 & 0.570 & 0.232 & 0.198  \\
   $\mu_1=5$ & 2.32 & 0.82 & 0.86 & 0.573 & 0.230 & 0.197    \\
   $\mu_1=10$ & 2.80 & 0.59 & 0.61 & 0.641 & 0.190& 0.169  \\
  $\mu_1=30$ & 3.43 & 0.28 & 0.29 & 0.814 & 0.097 & 0.089  \\  \hline
\end{tabular}}
\end{table}

\begin{table}[h]
 \caption{Numerical results for $\e[L_i]$, $\e[W_i]$ and $Cor(L_i,L_j)$ when $\lambda=4$, $\mu_2=5$, $\mu_3=5$, $p_1=0.9999$, $p_2=p_3=0.00005$, $q_1=q_2=q_3=1$}
  \label{SL3 mu2=5, mu3=5 p1=1, p2=p3=0 and q1=q2=q3=1}
    \centering\footnotesize{
\begin{tabular}{|l||c|c|c|c|c|c|c|c|c|c|}
  \hline
  $\mu_1$  & $\e[L_1]$ & $\e[L_2]$ &$\e[L_3]$& $\e[W_1]$ & $\e[W_2]$& $\e[W_3]$ & $Cor(L_1,L_2)$ & $Cor(L_1,L_3)$ & $Cor(L_2,L_3)$ & $GI$ \\ \hline
  $\mu_1=3.4$ & 46.31 & 46.08 & 46.08 & 22.60 & 47.24 & 47.24 &  0.9999 & 0.9999 &0.9999 & 0.0011 \\
   $\mu_1=3.5$& 14.46 & 14.23 & 14.23 & 6.99 & 14.74 & 14.74 &  0.9990 & 0.9990 &0.9990  &   0.0036\\
   $\mu_1=4$ & 3.42 & 3.20 & 3.20 & 1.58 & 3.48 & 3.48 &  0.9824 & 0.9824 &0.9860 & 0.0149  \\
   $\mu_1=4.5$ & 2.02 & 1.81 & 1.81 & 0.90 & 2.06 & 2.06 &  0.9525 & 0.9525 & 0.9628  & 0.0248 \\
  $\mu_1=5$ & 1.47 & 1.26 & 1.26 & 0.63 & 1.49 & 1.49 &  0.9158 & 0.9158 & 0.9348  &  0.0351  \\
  $\mu_1=10$ & 0.48 & 0.33 & 0.33 & 0.17 & 0.55 & 0.55 &  0.6079 & 0.6079 & 0.6973 & 0.0877 \\
  $\mu_1=30$ & 0.17 & 0.09 & 0.09 & 0.05 & 0.32 & 0.32 &  0.3774 & 0.3774 & 0.4718 & 0.1524 \\  \hline
\end{tabular}}
\end{table}

\begin{table}[h]
 \caption{Numerical results for $\lambda^i_{eff}$ and $\gamma_i$ when $\lambda=4$, $\mu_2=5$, $\mu_3=5$, $p_1=0.9999$, $p_2=p_3=0.00005$, $q_1=q_2=q_3=1$}
  \label{SL3 gamma mu2=5, mu3=5 p1=1, p2=p3=0 and q1=q2=q3=1}
    \centering\footnotesize{
\begin{tabular}{|l||c|c|c|c|c|c|c|}
  \hline
  $\mu_1$  & $\lambda^1_{eff}$ & $\lambda^2_{eff}$ & $\lambda^3_{eff}$ & $\gamma_1$ & $\gamma_2$ & $\gamma_3$ \\ \hline
  $\mu_1=3.4$ & 2.05 & 0.975 & 0.975 & 0.607 & 0.1965 & 0.1965  \\
   $\mu_1=3.5$& 2.07 & 0.965 & 0.965 & 0.605 & 0.1975 & 0.1975  \\
   $\mu_1=4$& 2.17 & 0.915 & 0.915 & 0.60 & 0.20 & 0.20  \\
   $\mu_1=4.5$& 2.25 & 0.875 & 0.875 & 0.60 & 0.20 & 0.20  \\
   $\mu_1=5$ & 2.32 & 0.84 & 0.84 & 0.601 & 0.1995 & 0.1995    \\
   $\mu_1=10$ & 2.81 & 0.595 & 0.595 & 0.665 & 0.1675 & 0.1675  \\
  $\mu_1=30$ & 3.45 & 0.275 & 0.275 & 0.83 & 0.085 & 0.085  \\  \hline
\end{tabular}}
\end{table}

\begin{table}[h]
 \caption{Numerical results for $\e[L_i]$, $\e[W_i]$ and $Cor(L_i,L_j)$ when $\lambda=4$, $\mu_1=4.5$, $\mu_2=3$, $\mu_3=5$, $p_1=p_2=p_3=\frac{1}{3}$, $q_2=q_3=1$}
  \label{SL3 q1 mu1=4.5 mu2=3, mu3=5 p1=p2=p3=1/3 and q2=q3=1}
    \centering\footnotesize{
\begin{tabular}{|l||c|c|c|c|c|c|c|c|c|c|}
  \hline
  $q_1$  & $\e[L_1]$ & $\e[L_2]$ &$\e[L_3]$& $\e[W_1]$ & $\e[W_2]$& $\e[W_3]$ & $Cor(L_1,L_2)$ & $Cor(L_1,L_3)$ & $Cor(L_2,L_3)$ & $GI$\\ \hline
  $q_1=1$ & 48.64 & 48.69 & 48.62 & 35.46 & 39.90 & 34.53 &  0.9999 & 0.9999 &0.9999 & 0.0003 \\
   $q_1=3$& 38.08 & 38.16 & 38.10 & 26.88 & 31.83 & 27.53 &  0.9998 & 0.9998 &0.9998  & 0.0005\\
   $q_1=8$ & 34.01 & 34.10 & 34.04 & 23.59 & 28.73 & 24.82 &  0.9998 & 0.9998 &0.9998  & 0.0006\\
   $q_1=13$ & 33.00 & 33.10 & 33.04 & 22.79 & 27.96 & 24.16 &  0.9998 & 0.9998 & 0.9998  & 0.0007 \\
  $q_1=20$ & 32.43 & 32.53 & 32.46 & 22.33 & 27.52 & 23.77 &  0.9998 & 0.9998 & 0.9998   &  0.0007\\
  $q_1=30$ & 32.07 & 32.17 & 32.10 & 22.04 & 27.25 & 23.54 &  0.9998 & 0.9998 & 0.9998  & 0.0007 \\
  $q_1=100$ & 31.56 & 31.66 & 31.60 & 21.63 & 26.86 & 23.20 &  0.9998 & 0.9998 & 0.9998  & 0.0007\\  \hline
\end{tabular}}
\end{table}

\begin{table}[h]
 \caption{Numerical results for $\lambda^i_{eff}$ and $\gamma_i$ when $\lambda=4$, $\mu_1=4.5$, $\mu_2=3$, $\mu_3=5$, $p_1=p_2=p_3=\frac{1}{3}$, $q_2=q_3=1$}
  \label{SL3 gamma q1 mu1=4.5 mu2=3, mu3=5 p1=p2=p3=1/3 and q2=q3=1}
    \centering\footnotesize{
\begin{tabular}{|l||c|c|c|c|c|c|c|}
  \hline
  $q_1$  & $\lambda^1_{eff}$ & $\lambda^2_{eff}$ & $\lambda^3_{eff}$ & $\gamma_1$ & $\gamma_2$ & $\gamma_3$ \\ \hline
  $q_1=1$ & 1.37 & 1.22 & 1.41 & 0.307 & 0.409 & 0.284  \\
   $q_1=3$& 1.42 & 1.20 & 1.38 & 0.318 & 0.402 & 0.280  \\
   $q_1=8$& 1.44 & 1.19 & 1.37 & 0.323 & 0.398 & 0.279  \\
   $q_1=13$& 1.45 & 1.18 & 1.37 & 0.323 & 0.398 & 0.279 \\
   $q_1=20$ & 1.45 & 1.18 & 1.37 & 0.326 & 0.397 & 0.277    \\
   $q_1=30$ & 1.46 & 1.18 & 1.36 & 0.326 & 0.397 & 0.277  \\
  $q_1=100$ & 1.46 & 1.18 & 1.36 & 0.327 & 0.396 & 0.277  \\  \hline
\end{tabular}}
\end{table}

\begin{table}[h]
 \caption{Numerical results for $\e[L_i]$, $\e[W_i]$ and $Cor(L_i,L_j)$ when $\lambda=4$, $\mu_1=4.2$, $\mu_2=3$, $\mu_3=5$, $p_1=0.9999$, $p_2=p_3=0.00005$, $q_2=q_3=1$}
  \label{SL3 q1 mu1=4.5 mu2=3, mu3=5 p1=1, p2=p3=0 and q2=q3=1}
    \centering\footnotesize{
\begin{tabular}{|l||c|c|c|c|c|c|c|c|c|c|}
  \hline
  $q_1$  & $\e[L_1]$ & $\e[L_2]$ &$\e[L_3]$& $\e[W_1]$ & $\e[W_2]$& $\e[W_3]$ & $Cor(L_1,L_2)$ & $Cor(L_1,L_3)$ & $Cor(L_2,L_3)$ & $GI$ \\ \hline
  $q_1=1$ & 105.10 & 104.95 & 104.90 & 47.79 & 123.78 & 110.07 &  0.9999 & 0.9999 &0.9999 & 0.0004\\
   $q_1=3$& 83.44 & 83.30 & 83.25 & 37.32 & 100.22 & 89.25 &  0.9999 & 0.9999 &0.9999 & 0.0005\\
   $q_1=8$ & 74.82 & 74.87 & 74.82 & 33.25 & 91.06 & 81.15 &  0.9999 & 0.9999 &0.9999 & 0.00001\\
   $q_1=13$ & 72.90 & 72.77 & 72.72 & 32.24 & 88.79 & 79.14 &  0.9999 & 0.9999 & 0.9999 & 0.0005 \\
  $q_1=20$ & 71.71 & 71.58 & 71.53 & 31.66 & 87.50 & 77.99 &  0.9999 & 0.9999 & 0.9999 & 0.0006\\
  $q_1=30$ & 70.96 & 70.83 & 70.78 & 31.30 & 86.68 & 77.27 &  0.9999 & 0.9999 & 0.9999  & 0.0006 \\
  $q_1=100$ & 69.89 & 69.78 & 69.72 & 30.79 & 85.53 & 76.26 &  0.9999 & 0.9999 & 0.9999 & 0.0005\\  \hline
\end{tabular}}
\end{table}

\begin{table}[h]
 \caption{Numerical results for $\lambda^i_{eff}$ and $\gamma_i$ when $\lambda=4$, $\mu_1=4.2$, $\mu_2=3$, $\mu_3=5$, $p_1=0.9999$, $p_2=p_3=0.00005$, $q_2=q_3=1$}
  \label{SL3 gamma q1 mu1=4.2 mu2=3, mu3=5 p1=1, p2=p3=0 and q2=q3=1}
    \centering\footnotesize{
\begin{tabular}{|l||c|c|c|c|c|c|c|}
  \hline
  $q_1$  & $\lambda^1_{eff}$ & $\lambda^2_{eff}$ & $\lambda^3_{eff}$ & $\gamma_1$ & $\gamma_2$ & $\gamma_3$ \\ \hline
  $q_1=1$ & 2.20 & 0.85 & 0.95 & 0.526 & 0.283 & 0.191  \\
   $q_1=3$& 2.24 & 0.83 & 0.93 & 0.535 & 0.278 & 0.187  \\
   $q_1=8$& 2.26 & 0.82 & 0.92 & 0.540 & 0.275 & 0.185  \\
   $q_1=13$& 2.26 & 0.82 & 0.92 & 0.541 & 0.274 & 0.185 \\
   $q_1=20$ & 2.26 & 0.82 & 0.92 & 0.542 & 0.274 & 0.184    \\
   $q_1=30$ & 2.27 & 0.82 & 0.91 & 0.543 & 0.273 & 0.184  \\
  $q_1=100$ & 2.27 & 0.82 & 0.91 & 0.543 & 0.273 & 0.184  \\  \hline
\end{tabular}}
\end{table}

\section{Conclusions}\label{SL3 sec conclusions}
This paper studies a queueing system combining two operating policies: Serve the Longest Queue and Join the Shortest Queue, that were separately studied by J.W. Cohen (\cite{cohen1987two} and \cite{cohen1998analysis}), respectively. The combined JSQ-SLQ model (first introduced in \cite{perel2020polling} for a preemptive 2-queue system) is analyzed for a 3-queue, single-server, Markovian system. The system is formulated in an un-conventional way, where, instead of defining an un-bounded 3-dimensional state space $(L_1,L_2,L_3)$, where $L_i$ represents the number of customers in $Q_i$, $i=1,2,3$, the system is characterized by the couple $L_1$ and $D=(I,L_1-L_2,L_1-L_3)$, where $I=1, 2, 3$ indicates the location of the server. This leads to a 2-dimensional state space with infinite dimension for $L_1$ but finite dimension for $D$. Consequently, one can apply the probability generating functions method, as well as the matrix geometric analysis, for QBD processes to derive the system's joint 2-dimensional probability mass function, the system's stability condition, and consequently, calculate its performance measures. Numerical results are presented in various tables with wide range of parameter values. It is shown that the main purpose of the JSQ-SLQ policy, namely balancing between the queue sizes, is greatly achieved. This is further demonstrated by the Gini Index, which is close to zero in most cases, while its highest value, $GI=0.1524$, is in the extreme case when $p_1=0.9999$ (Table \ref{SL3 mu2=5, mu3=5 p1=1, p2=p3=0 and q1=q2=q3=1}). Note that $GI=0$ in the symmetric case (Table \ref{SL3 mu2=5, mu3=5, p1=p2=p3=1/3 and q1=q2=q3=1}).

A natural extension of the presented model is to consider a preemptive single-server JSQ-SLQ system with $N\geq4$ non-symmetrical queues. Such a system can be modeled as a 2-dimensional continuous-time Markovian process, $(L_1,D)$, where $L_1$ denotes the (possibly infinite) size of $Q_1$, while $D$ is a vector where its first coordinate indicates the position of the server, and its $j$'th coordinate, $2\leq j\leq N$, represents the difference $L_1-L_j$. Although $D$ assumes values belonging to a finite set $\mathfrak{D}$, the size of $\mathfrak{D}$ increases exponentially with $N$. Specifically, consider the case when $L_1=0$. Then, either all queues are empty and the server may reside in any of the $N$ queues, or there are exactly $i$ queues with a single customer, $1\leq i\leq N-1$, where the server may reside in any of these $i$ queues. Therefore, the size of the set $\mathfrak{D}$ when $L_1=0$ is
\[
N+\sum_{i=1}^{N-1}{i\binom{N-1}{i}}=N+(N-1)\cdot2^{N-2}.
\]
For the case where $L_1\geq1$, the number of states in $\mathfrak{D}$ is $N\cdot2^{N-1}$. This follows since the server may reside in any of the $N$ queues, say $Q_i$, while, the number of customers in all other $N-1$ queues $Q_j$, $j\neq i$, can only be equal to $L_i$ or to $L_i-1$. Thus, the resulting system becomes analytically intractable.
\clearpage
\bibliographystyle{apalike}
\bibliography{shortest_longest_ref}
\end{document}